\documentclass[12pt,a4paper]{amsart}
\usepackage{fullpage,setspace}
\usepackage[T1]{fontenc}
\usepackage{lmodern}
\usepackage{amsmath,amssymb,amsthm,mathtools}
\usepackage{tikz,tikz-cd}
\usetikzlibrary{arrows.meta,calc,positioning,decorations.markings}
\usepackage{microtype}
\usepackage{xcolor}
\definecolor{dark-red}{rgb}{0.5,0.15,0.15}
\usepackage[colorlinks=true,linkcolor=black,citecolor=dark-red,urlcolor=dark-red,linktoc=all]{hyperref}
\usepackage[figure]{hypcap}
\AddToHook{env/equation/begin}{}
\hypersetup{pdftitle={Homotopical interpretation of directed Van Kampen},pdfauthor={Philippe Gaucher},pdfdisplaydoctitle=true}
\swapnumbers
\newtheorem{theorem}{Theorem}[section]
\newtheorem*{introVanKampen}{Homotopical interpretation of directed Van Kampen}
\newtheorem*{introModel}{The Grandis model structure}
\newtheorem*{grandisInvariance}{Grandis's invariance theorem}
\newtheorem*{smithRecognition}{Smith's recognition theorem}
\newtheorem{proposition}[theorem]{Proposition}
\newtheorem*{collarproposition}{Proposition}

\newtheorem{corollary}[theorem]{Corollary}
\theoremstyle{definition}
\newtheorem{remark}[theorem]{Remark}
\newtheorem{example}[theorem]{Example}
\renewcommand{\leq}{\leqslant}
\renewcommand{\geq}{\geqslant}
\renewcommand{\vec}{\overrightarrow}
\newcommand{\D}{\mathbf{dTop}}
\newcommand{\Cat}{\mathbf{Cat}}
\newcommand{\I}{\mathcal I_{\mathrm G}}
\newcommand{\MG}{\D_{\mathrm G}}
\newcommand{\WG}{W_{\mathrm G}}
\newcommand{\rev}[1]{#1_{\mathrm{nat}}}
\newcommand{\Iv}{\rev{I}}
\newcommand{\cube}[1]{C_{#1}}
\newcommand{\abs}[1]{\lvert#1\rvert}
\newcommand{\PiC}{\vec{\pi}_1}
\newcommand{\clp}[1]{\langle #1\rangle}
\DeclareMathOperator*{\hocolim}{hocolim}
\newcommand{\theoremlink}[1]{\hyperref[#1]{Theorem~\ref*{#1}}}

\newcommand{\GrandisModelStatement}{%
Let $\D$ denote either the category of directed spaces with $\Delta$-generated underlying spaces or its full subcategory of $\Delta$-Hausdorff objects. Then $\D$ admits a combinatorial model structure $\MG$ for which the seven maps in~\eqref{eq:generators} form a set of generating cofibrations. Its weak equivalences are precisely the directed maps $f$ for which $\PiC(f)$ is an equivalence of categories.
}
\newcommand{\VanKampenStatement}{%
The fundamental category functor $\PiC$ induces a functor from the underlying $\infty$-category of $\MG$ to that of the folk model category of small categories. This functor preserves all small colimits. In particular, for a diagram $X\colon J\to\MG$ indexed by a small category, the canonical comparison
\[
 \hocolim_{j\in J}\PiC(X(j))\longrightarrow\PiC\!\left(\hocolim_{j\in J}X(j)\right)
\]
is represented by an equivalence of categories, where $\hocolim$ denotes the homotopy colimit in the folk model structure on the left and in the Grandis model structure on the right.
}
\title{Homotopical interpretation of directed Van Kampen}
\author[P. Gaucher]{Philippe Gaucher}
\address{Universit\'e Paris Cit\'e, CNRS, IRIF, F-75013, Paris, France}
\urladdr{\url{https://www.irif.fr/~gaucher}}
\begin{document}

\begin{abstract}
We construct the Grandis model structure on directed spaces with $\Delta$-generated underlying spaces and on the full subcategory of $\Delta$-Hausdorff objects. In either setting the model structure is left proper, and its weak equivalences are precisely the maps that induce equivalences of fundamental categories in the sense of Grandis. The fundamental category functor preserves all small homotopy colimits with values in the folk model structure on small categories. For a union of two induced directed subspaces, the union square is a homotopy pushout precisely when Grandis's categorical pushout formula holds. His covering theorem gives this conclusion whenever the interiors of the two subspaces cover the space.
\end{abstract}

\maketitle
\setcounter{tocdepth}{1}
\tableofcontents
\hypersetup{linkcolor=dark-red}

\section{Introduction}

The purpose of this paper is to give a homotopical interpretation of Grandis's directed Van Kampen theorem. The fundamental category $\PiC(X)$ of a directed space $X$ has all points of $X$ as objects. Its arrows are $d$-homotopy classes of directed paths with fixed endpoints. A directed map $f\colon X\to Y$ induces the functor $\PiC(f)$ by sending points and path classes to their images under $f$. An elementary $d$-homotopy is a continuous map from the parameter square which fixes the endpoints and sends every coordinatewise nondecreasing path to a directed path. Two paths represent the same arrow if a finite zigzag of elementary $d$-homotopies joins them. A dihomotopy requires only a continuous family of directed paths with fixed endpoints~\cite[Definitions~4.34 and~4.49--4.50]{DAT_book}.

For a span $B\leftarrow A\to C$, write $B\sqcup_A C$ for its pushout. Grandis proved that if the interiors of two induced directed subspaces $U,V$ cover $X$, the canonical functor
\[
 \PiC(U)\sqcup_{\PiC(U\cap V)}\PiC(V)\longrightarrow\PiC(X)
\]
is an isomorphism~\cite[Theorem~3.6(a)]{mg} (Theorem~\ref{thm:grandis}). The inclusions of $U\cap V$ induce functors injective on objects. These are cofibrations in the folk model structure on small categories, in which every object is cofibrant. Thus the pushout on the left computes a homotopy pushout. In which model structure on directed spaces does the union $X=U\cup V$ itself compute the homotopy pushout of $U\leftarrow U\cap V\to V$?

We keep Grandis's definition of the fundamental category and choose the generating cofibrations with this question in mind. The directed interval represents paths, and the directed square represents the elementary homotopies relating paths when two opposite edges are constant. There is also a role for the ordinary interval and disk with all continuous paths directed. They allow us to compare different lifts of the same path by isomorphisms at their endpoints and to retract a collar of an attaching sphere. A directed map from the ordinary interval with all paths directed to the directed interval must be constant, since it must be nondecreasing in both orientations. This explains why the collar requires a second kind of parameter. These two kinds of parameter enter separately in the seven generators of~\eqref{eq:generators}.

\begin{introModel}[\theoremlink{thm:model}]
\GrandisModelStatement
\end{introModel}

For diagrams of ordinary spaces indexed by a small category, Dror Farjoun proved that the fundamental groupoid of the homotopy colimit is equivalent to the homotopy colimit of the fundamental groupoids, taken in groupoids~\cite[Theorem~1.1]{Farjoun2004}. Taking the automorphism group at a basepoint gives the classical fundamental group of the homotopy colimit. The following theorem gives a directed analogue for Grandis's fundamental category.

\begin{introVanKampen}[\theoremlink{thm:positive}]
\VanKampenStatement
\end{introVanKampen}

For an arbitrary union $X=U\cup V$ of induced directed subspaces, the square of inclusions is a homotopy pushout in the Grandis model structure precisely when Grandis's categorical comparison is an isomorphism. His theorem supplies this isomorphism when the interiors cover $X$. No cofibrancy assumption on $U$, $V$ or their intersection enters this interpretation (Corollary~\ref{cor:cover}). Two closed subintervals of the directed interval meeting at an endpoint still give a homotopy pushout (Example~\ref{ex:closed-halves}), whereas the alternating closed cover of Example~\ref{ex:alternating-cover} does not.

\subsection*{Outline of the paper}
Section~\ref{sec:setting} recalls the fundamental category and Grandis's covering theorem and introduces the generating cofibrations. Section~\ref{sec:lifts} examines the lifts of paths and their classes. In Section~\ref{sec:excision} we use the covering theorem to calculate the fundamental category after attaching a cell and extend the calculation to pushouts along cofibrations. These two calculations give the model structure in Section~\ref{sec:model}. Section~\ref{sec:van-kampen} proves the homotopy colimit theorem and applies it to unions.

\subsection*{Acknowledgment}
With the help of GPT Astra, I began by studying the minimal model structure generated by the boundary inclusions of the directed $n$-cubes, for $n\geq 0$. There is a trivial cofibration that identifies two distinct $d$-homotopy classes of paths with the same endpoints. Some partial fixes to this model structure and to the statement seemed possible, but none was satisfactory. I began to believe that the problem had no solution. Then, without much hope, I asked Astra whether we could change the generating cofibrations so that the weak equivalences would be precisely the maps inducing equivalences of fundamental categories in the sense of Grandis. To my surprise, the seven maps in~\eqref{eq:generators} did the trick (Theorem~\ref{thm:model}). The change is already visible for the directed interval (Remark~\ref{rem:interval-boundary}). I am grateful to Astra for helping me find the answer.

\section{Directed spaces and fundamental categories}\label{sec:setting}

Let $\mathbf{Top}$ denote either the category of $\Delta$-generated spaces or its full subcategory of $\Delta$-Hausdorff spaces. A space $K$ is \emph{$\Delta$-generated} if its topology is final for the continuous maps $\Delta^n\to K$ from the standard topological simplices, $n\geq0$. It is \emph{$\Delta$-Hausdorff} if the image of every continuous map $[0,1]\to K$ is closed~\cite[Definition~B.1]{leftproperflow}.

In either category, products and subspaces carry the $\Delta$-kelleyfication of the ordinary product or subspace topology. This is the final topology for the same continuous test simplices. Thus a continuous map from a $\Delta$-generated space whose image lies in a subspace factors continuously through that subspace. We use this in particular for paths, squares and simplices. Products and subspaces of $\Delta$-Hausdorff spaces are again $\Delta$-Hausdorff~\cite[Propositions~B.4 and~B.5]{leftproperflow}.

Let $\D$ denote the category of directed spaces over the chosen $\mathbf{Top}$. A \textit{directed space} $X$ consists of an underlying space $\abs X$ in $\mathbf{Top}$ and a set of continuous paths $[0,1]\to\abs X$, called its directed paths. This set contains the constant paths and is closed under concatenation and precomposition with any continuous nondecreasing map $[0,1]\to[0,1]$. A directed map is a continuous map preserving directed paths. In Grandis's terminology, the \emph{natural directed structure} on a space $K$ is the one in which every continuous path $[0,1]\to K$ is directed~\cite[Section~1.1, pp.~284--285]{mg}. We write $\rev K$ for $K$ with this structure. The directed paths of an induced subspace are the directed paths of the ambient space whose images lie in the subspace. On a product, a path is directed when both projections are directed. For either choice of $\mathbf{Top}$, the category $\D$ is locally presentable~\cite[p.~21, after Definition~3.5]{GlobularNaturalSystem}.

To form a colimit in $\D$, take the colimit of the underlying spaces in $\mathbf{Top}$ and give it the directed structure generated by the images of the directed paths in the diagram. Colimits of $\Delta$-generated spaces have the ordinary final topology. In the $\Delta$-Hausdorff setting we then take the $\Delta$-Hausdorff reflection, whose unit is a quotient map~\cite[Proposition~B.7 and its proof]{leftproperflow}. Thus the images of the summands still cover the underlying set. In either setting a path in a pushout is directed precisely when it admits a finite subdivision on which each piece is the image of a directed path in one summand. Indeed, paths with this property include the images of directed paths from the summands and the constant paths. They are closed under concatenation and under nondecreasing reparametrization, since the inverse images of the subdivision points under a nondecreasing map are intervals or empty. Conversely, every such path is a reparametrized concatenation of images of directed paths.

For $n\geq0$ let $\cube n$ be the unit cube $[0,1]^n$ with its coordinatewise nondecreasing continuous paths, where $\cube 0=*$ is the one-point space. For $I=[0,1]$, Grandis calls $\Iv$ the \emph{natural interval}~\cite[Section~1.2]{mg}. We use $\cube 1$ for the directed path parameter and $\Iv$ for homotopies whose parameter can run in either direction.

\begin{figure}[htbp]
\centering
\begin{tikzpicture}[x=1cm,y=1cm,>=Stealth,font=\small,
  oriented/.style={postaction={decorate},decoration={markings,mark=at position .55 with {\arrow{Stealth[length=1.9mm]}}}},
  side/.style={thick,oriented},slice/.style={black!65,dashed,oriented},
  track/.style={dark-red,oriented,thick}]
  \begin{scope}
    \fill[black!3] (0,0) rectangle (2.6,2.6);
    \draw[side] (0,0) -- (2.6,0) node[midway,below=3pt] {$t=0$};
    \draw[side] (0,2.6) -- (2.6,2.6) node[midway,above=3pt] {$t=1$};
    \draw[side] (0,0) -- (0,2.6) node[midway,left=3pt] {$s=0$};
    \draw[side] (2.6,0) -- (2.6,2.6) node[midway,right=3pt] {$s=1$};
    \draw[slice] (1.3,0) -- (1.3,2.6);
    \draw[track] (0,1.04) -- (2.6,1.04);
  \end{scope}
  \draw[->,semithick] (4,1.3) -- (5.15,1.3) node[midway,above] {$H$};
  \begin{scope}[shift={(7.25,0)}]
    \draw[side] (0,0) .. controls (-1.52,.72) and (-1.52,1.88) .. (0,2.6)
      coordinate[pos=.4] (lefttrack) node[pos=.53,left=3pt] {$\gamma_0$};
    \draw[side] (0,0) .. controls (1.52,.72) and (1.52,1.88) .. (0,2.6)
      coordinate[pos=.4] (righttrack) node[pos=.53,right=3pt] {$\gamma_1$};
    \draw[slice] (0,0) -- (0,2.6);
    \draw[track] (lefttrack) -- (righttrack);
    \fill (0,0) circle (1.7pt) node[below=3pt] {$x$};
    \fill (0,2.6) circle (1.7pt) node[above=3pt] {$y$};
  \end{scope}
\end{tikzpicture}
\caption{An elementary $d$-homotopy, with its image in $X$ drawn schematically. The dashed segment in the parameter square gives an intermediate path from $x$ to $y$. The horizontal arrow inside the square shows how a point on the path moves as $s$ increases with $t$ fixed. Every coordinatewise nondecreasing path in the square must have directed image in $X$.}\label{fig:elementary-d-homotopy}
\end{figure}
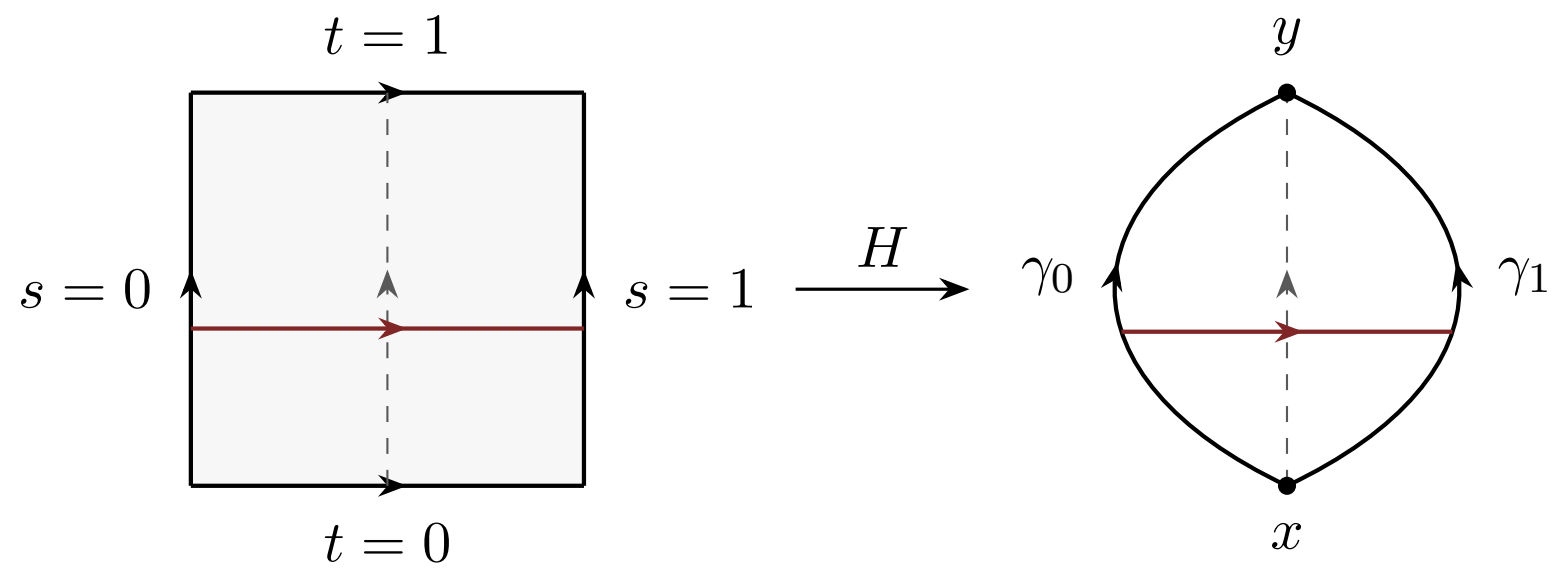

For objects $x,y$ of a category $\mathcal C$, write $\mathcal C(x,y)$ for the set of arrows from $x$ to $y$ and $\mathrm{id}_x$ for the identity of $x$. The objects of the fundamental category $\PiC(X)$ are the points of $X$. Write $\clp{\gamma}$ for the arrow represented by a directed path $\gamma$. Two paths $\gamma_0,\gamma_1$ from $x$ to $y$ represent the same arrow if and only if they are joined by a finite zigzag of elementary $d$-homotopies. Such a homotopy from $\gamma_0$ to $\gamma_1$ (Figure~\ref{fig:elementary-d-homotopy}) is a directed map $H\colon\cube 2\to X$ with
\[
 H(0,t)=\gamma_0(t),\quad H(1,t)=\gamma_1(t),\quad H(s,0)=x,\quad H(s,1)=y \qquad (s,t\in[0,1]).
\]

\begin{grandisInvariance}[{\cite[Theorem~3.2(a),(c),(d)]{mg}}]
With composition induced by concatenation and identities represented by constant paths, the fundamental category defines a functor $\PiC\colon\D\to\Cat$ to the category of small categories. For a directed map $f\colon X\to Y$, the functor $\PiC(f)$ sends $x$ to $f(x)$ and $\clp{\gamma}$ to $\clp{f\circ\gamma}$. The functor $\PiC$ preserves small coproducts and binary products.

Let $H\colon X\times\cube1\to Y$ be a directed map, and define $H_\varepsilon\colon X\to Y$ by $H_\varepsilon(x)=H(x,\varepsilon)$ for $\varepsilon=0,1$. This $d$-homotopy induces a natural transformation $\eta\colon\PiC(H_0)\Rightarrow\PiC(H_1)$ between functors $\PiC(X)\to\PiC(Y)$. If $H$ is directed on $X\times\Iv$, then $\eta$ is a natural isomorphism. For an arrow $u\colon x\to y$ in $\PiC(X)$, the components satisfy
\begin{equation}\label{eq:homotopy}
 \eta_x=\clp{t\mapsto H(x,t)}\colon H_0(x)\to H_1(x),\qquad
 \eta_y\circ\PiC(H_0)(u)=\PiC(H_1)(u)\circ\eta_x.
\end{equation}
\end{grandisInvariance}

For two directed paths $\gamma_0,\gamma_1$ in $\cube1$ with the same endpoints, let $m(t)=\max\{\gamma_0(t),\gamma_1(t)\}$. For $i=0,1$, the formula
\[
 H_i(s,t)=(1-s)\gamma_i(t)+s\,m(t),\qquad 0\leq s,t\leq1,
\]
gives a $d$-homotopy from $\gamma_i$ to $m$ with fixed endpoints. The expression is nondecreasing in $s$ since $\gamma_i(t)\leq m(t)$, and in $t$ since both paths are nondecreasing. Taking the identity and a nondecreasing reparametrization fixing $0,1$ as paths in $\cube 1$, and composing these homotopies with a directed path in $X$, shows that such a reparametrization preserves its class. Since $\cube n$ is the $n$-fold product of $\cube1$ and $\cube0$ is a point, preservation of binary products gives
\begin{equation}\label{eq:thin}
 \PiC(\cube n)\cong\bigl([0,1]^n,\leq\bigr).
\end{equation}
On the right we regard the coordinatewise order as a category, with a unique arrow $x\to y$ when $x\leq y$ and none otherwise. A directed path can join $x$ to $y$ only if $x\leq y$, and then the straight segment is such a path~\cite[Section~3.4(a),(c)]{mg}.

For $a,b\in I$, the straight path $t\mapsto(1-t)a+tb$ represents an arrow $a\to b$ of $\PiC(\Iv)$. This arrow is unique because straight-line interpolation gives a homotopy between any two paths with these endpoints.

For an induced directed subspace $U$ of $X$, let $\operatorname{int}_X U$ denote the largest open subset of $\abs X$ contained in $\abs U$.

\begin{theorem}[Grandis's directed Van Kampen theorem {\cite[Theorem~3.6(a)]{mg}}]\label{thm:grandis}
Let $U,V$ be induced directed subspaces of $X$ with $\operatorname{int}_X U\cup\operatorname{int}_X V=\abs X$. Put $W=U\cap V$. The following square is a pushout of small categories.
\[
\begin{tikzcd}
 \PiC(W)\ar[r]\ar[d]&\PiC(U)\ar[d]\\
 \PiC(V)\ar[r]&\PiC(X)\ar[ul,phantom,"\ulcorner",very near start]
\end{tikzcd}
\]
\end{theorem}

In the folk model structure on the category $\Cat$ of small categories, the weak equivalences are the equivalences of categories, the cofibrations are the functors injective on objects and the fibrations are the isofibrations~\cite[Examples~11.3.9 and~12.2.1]{Riehl}. In particular all small categories are cofibrant. The folk model structure is combinatorial since it is cofibrantly generated and $\Cat$, the category of models of the finite-limit theory of categories, is locally finitely presentable~\cite[Theorem~5.30]{TheBook}. Write $\varnothing$ for the initial object. For a span $B\leftarrow A\to C$, the latching maps in the direct Reedy model structure are $\varnothing\to A$, $A\to B$ and $A\to C$. If they are cofibrations, the span is cofibrant and its ordinary pushout computes the homotopy pushout~\cite[Theorem~5.1.3 and Corollary~5.1.6]{MR99h:55031}. This applies to Theorem~\ref{thm:grandis}, since both inclusions of $W$ are injective on points.

\subsection*{The generating cofibrations}
For each integer $m\geq1$, let $\mathsf{D}^m$ be the closed unit disk in $\mathbb R^m$ and let $\mathsf{S}^{m-1}$ be its boundary. We identify $\mathsf{D}^1$ with $[0,1]$ by the increasing affine homeomorphism and write $\mathsf{S}^0$ for its two endpoints with their discrete topology and hence their unique directed structure. The set $\I$ of generating cofibrations consists of the following maps in $\D$.
\begin{equation}\label{eq:generators}
\begin{gathered}
 \varnothing\longrightarrow\cube 0,\qquad
 \varnothing\longrightarrow\cube 1,\qquad
 \varnothing\longrightarrow\cube 2,\qquad
 \varnothing\longrightarrow\rev{\mathsf{S}^1},\\[2pt]
 \mathsf{S}^0\hookrightarrow\Iv,\qquad
 \cube 1\times \mathsf{S}^0\hookrightarrow\cube 1\times\Iv,\qquad
 \rev{\mathsf{S}^1}\hookrightarrow\rev{\mathsf{D}^2}.
\end{gathered}
\end{equation}
A relative $\I$-cell map is a transfinite composite of pushouts of coproducts of maps in $\I$. We write $\operatorname{cof}(\I)$ for the closure of $\I$ under pushouts, transfinite compositions and retracts. We call its elements $\I$-cofibrations. An object $K$ is cofibrant for these generators when $\varnothing\to K$ belongs to $\operatorname{cof}(\I)$. Write $\operatorname{inj}(\I)$ for the class of maps with the right lifting property against every member of $\I$, and $\WG$ for the class of directed maps inducing equivalences under $\PiC$.

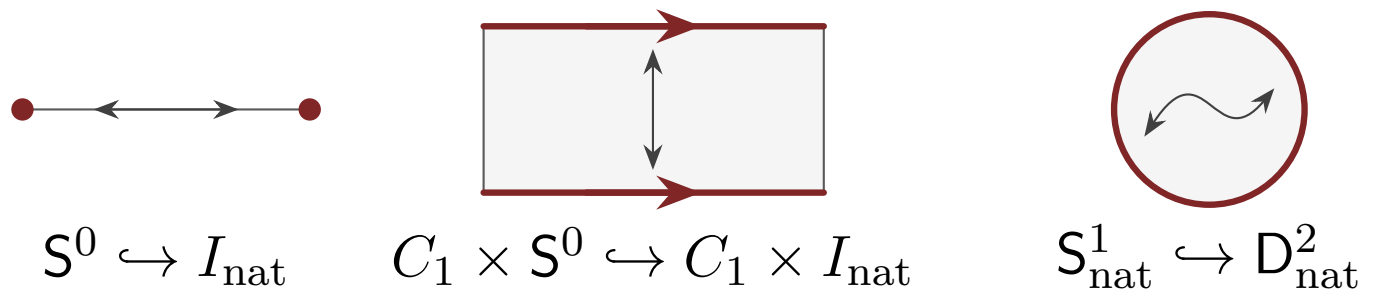
\begin{figure}[htbp]
\centering
\begin{tikzpicture}[>=Stealth,line cap=round,line join=round,font=\small]
 \tikzset{boundary/.style={draw=dark-red,very thick},direction/.style={draw=black!75,thin}}
 \begin{scope}
  \draw[black!65] (0,0)--(1.9,0);
  \draw[direction,<->] (.47,0)--(1.43,0);
  \foreach \x in {0,1.9}{\fill[dark-red] (\x,0) circle (2.1pt);}
  \node at (.95,-.97) {$\mathsf{S}^0\hookrightarrow\Iv$};
 \end{scope}
 \begin{scope}[shift={(3.05,0)}]
  \fill[black!4] (0,-.55) rectangle (2.25,.55);
  \draw[black!65] (0,-.55) rectangle (2.25,.55);
  \foreach \y in {-.55,.55}{
    \draw[boundary] (0,\y)--(2.25,\y);
    \draw[boundary,->] (.68,\y)--(1.44,\y);
  }
  \draw[direction,<->] (1.12,-.4)--(1.12,.4);
  \node at (1.12,-.97) {$\cube1\times \mathsf{S}^0\hookrightarrow\cube1\times\Iv$};
 \end{scope}
 \begin{scope}[shift={(7.85,0)}]
  \fill[black!4] (0,0) circle (.63);
  \draw[boundary] (0,0) circle (.63);
  \draw[direction,<->] (-.43,-.18) .. controls (-.03,.43) and (.04,-.33) .. (.43,.14);
  \node at (0,-.97) {$\rev{\mathsf{S}^1}\hookrightarrow\rev{\mathsf{D}^2}$};
 \end{scope}
\end{tikzpicture}
\caption{The three generators with nonempty source. The attaching subspace is drawn in red. In the strip a path may move freely in the vertical coordinate and must be nondecreasing horizontally. In the disk every continuous path is directed.}\label{fig:generators}
\end{figure}

Lifting against the first three maps amounts to lifting points, paths and directed squares, respectively, with no boundary data, as in the left square of~\eqref{eq:cube-strip-lifts}. The second line of~\eqref{eq:generators} supplies the comparisons illustrated in Figure~\ref{fig:generators}. The interval joins two lifts in one fibre, the strip compares lifts of a path, and the disk identifies the classes of two directed maps from $\Iv$ into one fibre whenever their endpoints agree. The generator $\varnothing\to\rev{\mathsf{S}^1}$ makes the source of the disk inclusion cofibrant with its full natural directed structure.

\section{Lifting paths and their classes}\label{sec:lifts}

A directed map $u\colon\Iv\to X$ gives a directed path by precomposition with the identity on points $\cube1\to\Iv$. We denote the class of this composite by $\clp u$. It is invertible, since the arrow $0\to1$ in $\PiC(\Iv)$ is invertible. The inverse is represented by $t\mapsto u(1-t)$. We shall use maps defined on $\Iv$ to join points of the same fibre. Grandis calls a path \emph{reversible} when its reverse is directed. Being directed on $\Iv$ is a stronger condition, and it need not be preserved by concatenation~\cite[Section~1.5]{mg}. We write $vu$ for an arrow $u$ followed by an arrow $v$.

\begin{proposition}\label{prop:fibre-transport}
Let $p\colon X\to Y$ be a directed map with the right lifting property against $\mathsf{S}^0\hookrightarrow\Iv$ and $\rev{\mathsf{S}^1}\hookrightarrow\rev{\mathsf{D}^2}$. For $y\in Y$, give the fibre $p^{-1}(y)=\{x\in X\mid p(x)=y\}$ its induced directed structure. For $a,b\in p^{-1}(y)$, there is a directed map $\Iv\to p^{-1}(y)$ from $a$ to $b$. All such maps represent the same isomorphism $\theta_{a,b}\colon a\to b$ in $\PiC(X)$. For points $a,b,c$ in one fibre,
\[
 \theta_{a,a}=\mathrm{id}_a,\qquad \theta_{b,c}\,\theta_{a,b}=\theta_{a,c}.
\]
\end{proposition}
\begin{proof}
In the following square the top map takes the endpoints of $\Iv$ to $a$ and $b$, and the bottom map is constant at $y$. The right lifting property of $p$ gives the indicated map $u$, which lands in the fibre.
\[
\begin{tikzcd}
 \mathsf{S}^0\ar[r,"{(a,b)}"]\ar[d,hook]&X\ar[d,"p"]\\
 \Iv\ar[r,"y"']\ar[ur,dashed,"u"]&Y .
\end{tikzcd}
\]

Consider directed maps $u_j\colon\Iv\to p^{-1}(y)$, $1\leq j\leq r$, whose underlying paths concatenate to a loop at $a$. We show that $\clp{u_r}\cdots\clp{u_1}=\mathrm{id}_a$. Insert a constant path if necessary so that $r\geq2$. Divide $\mathsf{S}^1$ into consecutive closed arcs $A_1,B_1,\ldots,A_r,B_r$ in cyclic order, starting at $z_0$. On $A_j$, traverse $u_j$ using an orientation-preserving parametrization by $[0,1]$. On $B_j$, take the constant value $u_j(1)=u_{j+1}(0)$, with indices read cyclically. These assignments agree at the endpoints and give a continuous map $L\colon\mathsf{S}^1\to\abs{p^{-1}(y)}$.

Give each arc $U_j=B_{j-1}\cup A_j\cup B_j$ its natural directed structure.\footnote{For arcs glued only at their endpoints, the directed colimit contains just the paths that subdivide finitely among the arcs. In an angular coordinate centered at a junction, the path $t\mapsto t\sin(1/t)$, extended by $0$ at $t=0$, visits both sides on every sufficiently short initial interval. It admits no such subdivision. The overlap puts a neighborhood of the junction in one arc. See~\cite[Section~1.5, pp.~289--290]{mg} for the interval example.} On $U_j$, the map $L$ factors through $u_j$ by collapsing $B_{j-1}$ and $B_j$ to the endpoints of the parameter interval of $A_j$. This collapse is continuous and hence directed between the natural structures. The interiors of the $U_j$ cover the circle, so every continuous path admits a finite subdivision into pieces lying in the $U_j$. Thus $\rev{\mathsf{S}^1}$ is the colimit in $\D$ of these arcs along their intersections. The restrictions of $L$ therefore give a directed map $\rev{\mathsf{S}^1}\to p^{-1}(y)$. Composing with the inclusion of the fibre gives the top map below. The right lifting property of $p$ against the disk inclusion gives the indicated extension $F$.
\[
\begin{tikzcd}
 \rev{\mathsf{S}^1}\ar[r,"{\mathrm{incl}\circ L}"]\ar[d,hook]&X\ar[d,"p"]\\
 \rev{\mathsf{D}^2}\ar[r,"y"']\ar[ur,dashed,"F"]&Y .
\end{tikzcd}
\]

Let $z\colon[0,1]\to\mathsf{S}^1$ traverse the arcs once in cyclic order from $z_0$ back to $z_0$. Since $L$ is constant on each $B_j$, the path $Lz$ represents $\clp{u_r}\cdots\clp{u_1}$. Convexity of the disk and its natural directed structure give a directed map
\[
 H\colon\cube2\longrightarrow X,\qquad H(s,t)=F\bigl((1-s)z(t)+s z_0\bigr).
\]
Indeed the argument of $F$ defines a continuous map into $\rev{\mathsf{D}^2}$, where every continuous path is directed. We have $H(0,t)=Lz(t)$ and $H(1,t)=H(s,0)=H(s,1)=a$. Thus the loop represents the identity. For two maps $u,v\colon\Iv\to p^{-1}(y)$ from $a$ to $b$, the loop obtained by following $u$ and then the reverse of $v$ gives $\clp v^{-1}\clp u=\mathrm{id}_a$. This proves uniqueness. The constant path represents $\theta_{a,a}$. Following paths from $a$ to $b$ and from $b$ to $c$, then the reverse of one from $a$ to $c$, gives $\theta_{a,c}^{-1}\theta_{b,c}\theta_{a,b}=\mathrm{id}_a$.
\end{proof}

\begin{proposition}\label{prop:lifts}
If $p\colon X\to Y$ has the right lifting property against $\I$, then $\PiC(p)$ is surjective on objects and fully faithful.
\end{proposition}
\begin{proof}
For $n=0,1,2$, lifting against $\varnothing\to\cube n$ in~\eqref{eq:generators} gives the map $\widetilde u$ on the left below. Thus $p$ is surjective on points, and every directed path or directed square in $Y$ has a lift. Given two lifts $\alpha,\beta\colon\cube1\to X$ of a directed path $\gamma\colon\cube1\to Y$, lifting against the strip inclusion gives $H$ on the right, where $\operatorname{pr}_1$ is projection onto $\cube1$.
\begin{equation}\label{eq:cube-strip-lifts}
\begin{tikzcd}
 \varnothing\ar[r]\ar[d,hook]&X\ar[d,"p"]\\
 \cube{n}\ar[r,"u"']\ar[ur,dashed,"\widetilde u"]&Y
\end{tikzcd}
\qquad
\begin{tikzcd}
 \cube{1}\times \mathsf{S}^0\ar[r,"{(\alpha,\beta)}"]\ar[d,hook]&X\ar[d,"p"]\\
 \cube{1}\times\Iv\ar[r,"\gamma\operatorname{pr}_1"']\ar[ur,dashed,"H"]&Y .
\end{tikzcd}
\end{equation}
By Grandis's invariance theorem, $H$ induces a natural isomorphism $\PiC(\alpha)\cong\PiC(\beta)$. By~\eqref{eq:homotopy}, its component at $i\in\{0,1\}$ is represented by the directed map $\Iv\to p^{-1}(\gamma(i))$ given by $s\mapsto H(i,s)$, so Proposition~\ref{prop:fibre-transport} identifies it with $\theta_{\alpha(i),\beta(i)}$. Naturality at $0\to1$ in $\PiC(\cube1)$ gives
\begin{equation}\label{eq:lift-comparison}
\begin{tikzcd}[column sep=large,row sep=large]
 \alpha(0)\ar[r,"\clp{\alpha}"]\ar[d,"{\theta_{\alpha(0),\beta(0)}}"']&
 \alpha(1)\ar[d,"{\theta_{\alpha(1),\beta(1)}}"]\\
 \beta(0)\ar[r,"\clp{\beta}"']&\beta(1).
\end{tikzcd}
\end{equation}
In particular, if $\delta\colon\cube1\to X$ lies in one fibre, take $\alpha$ constant at $\delta(0)$ and $\beta=\delta$. The top and left arrows are identities, so $\clp{\delta}=\theta_{\delta(0),\delta(1)}$. Here $\delta$ need only be directed on $\cube1$, as will be the case for the edges of a lifted square.

Fix $a,b\in X$. For a directed path $\gamma$ from $p(a)$ to $p(b)$ choose a lift $\alpha$ and put
\[
 T_{a,b}(\gamma)=\theta_{\alpha(1),b}\,\clp{\alpha}\,\theta_{a,\alpha(0)}
 \quad\text{in }\PiC(X)(a,b).
\]
For a second lift $\beta$, the composition law of Proposition~\ref{prop:fibre-transport} and the square~\eqref{eq:lift-comparison} give
\begin{equation}\label{eq:transport-calculation}
\begin{aligned}
 \theta_{\beta(1),b}\,\clp{\beta}\,\theta_{a,\beta(0)}
 &=\theta_{\beta(1),b}\,\clp{\beta}\,\theta_{\alpha(0),\beta(0)}\,\theta_{a,\alpha(0)}\\
 &=\theta_{\beta(1),b}\,\theta_{\alpha(1),\beta(1)}\,\clp{\alpha}\,\theta_{a,\alpha(0)}\\
 &=\theta_{\alpha(1),b}\,\clp{\alpha}\,\theta_{a,\alpha(0)}.
\end{aligned}
\end{equation}
Thus $T_{a,b}(\gamma)$ is independent of the chosen lift.

Now let $K\colon\cube2\to Y$ be an elementary $d$-homotopy from $\gamma_0$ to $\gamma_1$ fixing $p(a)$ and $p(b)$. The left square in~\eqref{eq:cube-strip-lifts} with $n=2$ and $u=K$ gives a lift $\widetilde K\colon\cube2\to X$. For $i=0,1$ the path $\alpha_i\colon\cube1\to X$ given by $\alpha_i(t)=\widetilde K(i,t)$ lifts $\gamma_i$. For $j=0,1$, the path $s\mapsto\widetilde K(s,j)$ lies in the fibre over $\gamma_0(j)=\gamma_1(j)$ and therefore represents $\theta_{\alpha_0(j),\alpha_1(j)}$. By~\eqref{eq:thin}, the two boundary routes from $(0,0)$ to $(1,1)$ in $\cube2$ have the same class. Applying $\PiC(\widetilde K)$ gives
\[
 \theta_{\alpha_0(1),\alpha_1(1)}\,\clp{\alpha_0}
 =\clp{\alpha_1}\,\theta_{\alpha_0(0),\alpha_1(0)}.
\]
Compose on the left with $\theta_{\alpha_1(1),b}$ and on the right with $\theta_{a,\alpha_0(0)}$. The composition law of Proposition~\ref{prop:fibre-transport} gives $T_{a,b}(\gamma_0)=T_{a,b}(\gamma_1)$ as in~\eqref{eq:transport-calculation}. The equality passes along every finite zigzag of elementary $d$-homotopies.

We can therefore put $\overline T_{a,b}(\clp{\gamma})=T_{a,b}(\gamma)$. The resulting maps
\[
\begin{tikzcd}[column sep=large]
 \PiC(X)(a,b) \ar[r,shift left=0.6ex,"\PiC(p)"] &
 \PiC(Y)(p(a),p(b)) \ar[l,shift left=0.6ex,"{\overline T_{a,b}}"]
\end{tikzcd}
\]
are inverse. Indeed, $p$ sends both fibre arrows in $T_{a,b}(\gamma)$ to identities, giving $\PiC(p)(T_{a,b}(\gamma))=\clp{\gamma}$. For a directed path $\delta$ from $a$ to $b$, taking $\delta$ itself as the lift of $p\delta$ gives $T_{a,b}(p\delta)=\clp{\delta}$. This proves full faithfulness.
\end{proof}

\section{Attaching a disk with its natural directed structure}\label{sec:excision}

Up to the identification of $\mathsf{D}^1$ with $I$, the three generators in~\eqref{eq:generators} with nonempty domain are the boundary inclusions $B\hookrightarrow E$, where
\[
 B=T\times\rev{\mathsf{S}^{m-1}},\qquad E=T\times\rev{\mathsf{D}^m},\qquad (T,m)\in\{(*,1),(\cube{1},1),(*,2)\}.
\]
Write $(t,z)$ for a point of $\abs T\times\mathsf{D}^m$ and $\|z\|$ for the Euclidean norm of $z$. In dimension one the interval coordinate $u\in[0,1]$ corresponds to $z=2u-1$, so both endpoints have radius $1$. For the strip $\cube{1}\times\rev{\mathsf{D}^1}$ the deformations below change only $z$ and leave the directed path coordinate fixed.

\begin{theorem}\label{thm:disk-excision}
Let $f\colon B\to A$ be a directed map and form the pushout on the left in $\D$. Then the square on the right is a pushout of small categories.
\[
\begin{tikzcd}
 B\ar[r,hook]\ar[d,"f"']&E\ar[d,"\chi"]\\
 A\ar[r,hook,"j"']&X\ar[ul,phantom,"\ulcorner",very near start]
\end{tikzcd}
\qquad
\begin{tikzcd}
 \PiC(B)\ar[r]\ar[d,"\PiC(f)"']&\PiC(E)\ar[d,"\PiC(\chi)"]\\
 \PiC(A)\ar[r,"\PiC(j)"']&\PiC(X)\ar[ul,phantom,"\ulcorner",very near start].
\end{tikzcd}
\]
\end{theorem}
The boundary $\abs B$ is a closed subspace of the Hausdorff space $\abs E$. In $\Delta$-generated spaces the underlying pushout therefore contains $\abs A$ as a closed subspace and is $\Delta$-Hausdorff whenever $\abs A$ is $\Delta$-Hausdorff~\cite[Proposition~B.15]{leftproperflow}. Thus in either setting $\abs X$ is the quotient of $\abs A\sqcup\abs E$ obtained by identifying each boundary point $(t,z)$ with $f(t,z)$. We identify $A$ with its image under $j$.

The continuous function $\rho\colon\abs X\to[0,1]$ which is $1$ on $A$ and takes $\chi(t,z)$ to $\|z\|$ is well defined because these values agree on $B$. Choose $0<a<b<1$ and give $U,W,V$ in~\eqref{eq:collar-spans} their induced directed structures. Here $i\colon B\hookrightarrow E$ is the boundary inclusion, and $\Phi\colon\mathsf{D}^m\to\mathsf{D}^m$ agrees with radial projection for $\|z\|\geq a$.
\begin{equation}\label{eq:collar-spans}
\begin{tikzcd}[column sep=1.65em,row sep=.22em]
 U=\rho^{-1}\bigl({]}a,1]\bigr)\ar[d,"d_U"'] &
 W=U\cap V\ar[l,hook]\ar[r,hook]\ar[d,"d_W"] &
 V=\rho^{-1}\bigl([0,b[\bigr)\ar[d,"d_V"] \\[1.25em]
 A & B\ar[l,"f"']\ar[r,hook,"i"] & E \\[.4em]
 {\chi(t,z)\longmapsto f\!\left(t,\dfrac{z}{\|z\|}\right)} &
 {\chi(t,z)\longmapsto \left(t,\dfrac{z}{\|z\|}\right)} &
 {\chi(t,z)\longmapsto(t,\Phi(z))} \\
 {x\longmapsto x\quad(x\in A)} & {} &
 {\Phi(z)=\dfrac{z}{\max(a,\|z\|)}}
\end{tikzcd}
\end{equation}
On $W$ we have $a<\|z\|<b$, so $\Phi(z)=z/\|z\|$. Thus $d_U|_W=f\circ d_W$ and $d_V|_W=i\circ d_W$. The map $\Phi$ expands the disk of radius $a$ to the unit disk (Figure~\ref{fig:reversible-collar}).

\begin{collarproposition}[Collar deformations]
The vertical maps in~\eqref{eq:collar-spans} are directed. For $Z=U,W,V$, the map $d_Z$ has a directed section $s_Z$ and there is a directed homotopy $H_Z\colon Z\times\Iv\to Z$ from $\mathrm{id}_Z$ to $s_Zd_Z$. Moreover, there is a directed homotopy $H_X\colon X\times\Iv\to X$ from $\mathrm{id}_X$ to a map $h$ such that $h|_U=jd_U$ and $h|_V=\chi d_V$.
\end{collarproposition}
\begin{proof}
\emph{Sections.} Choose $a<c<b$. On points, define
\[
\begin{alignedat}{2}
 s_U&\colon A\longrightarrow U, &\qquad s_U(x)&=j(x),\\
 s_W&\colon B\longrightarrow W, & s_W(t,z)&=\chi(t,cz),\\
 s_V&\colon E\longrightarrow V, & s_V(t,z)&=\chi(t,az).
\end{alignedat}
\]
For $\|z\|=1$ we have $\|cz\|=c\in{]}a,b[$ and $cz/\|cz\|=z$. For $\|z\|\leq1$ we have $\|az\|\leq a<b$ and $\Phi(az)=az/a=z$. Thus the sections take values in the indicated sets and~\eqref{eq:collar-spans} gives
\[
 d_U s_U=\mathrm{id}_A,\qquad d_W s_W=\mathrm{id}_B,\qquad d_V s_V=\mathrm{id}_E .
\]

\emph{Continuity and directedness.} The open sets $U$ and $V$ cover $X$. The attaching quotient $p=[j,\chi]\colon\abs A\sqcup\abs E\twoheadrightarrow\abs X$ restricts over the open set $\abs U$ to a quotient map
\[
 p_U\colon\abs A\sqcup\chi^{-1}(\abs U)\twoheadrightarrow\abs U,\qquad
 \chi^{-1}(\abs U)=\{(t,z)\in\abs E\mid a<\|z\|\leq1\}.
\]
Both summands are $\Delta$-generated, so their ordinary quotient already has the $\Delta$-generated topology of $U$. On the two summands one has
\[
 d_U p_U(x)=x,\qquad d_U p_U(t,z)=f\!\left(t,\frac{z}{\|z\|}\right).
\]
These continuous formulas agree at the attaching points since $z/\|z\|=z$ for $\|z\|=1$. Hence $d_U$ is continuous by the quotient property of $p_U$. By the description of pushout paths in Section~\ref{sec:setting}, a directed path in $U$ subdivides into images of directed paths in $A$ or the collar. On a piece from $A$, $d_U$ is the identity. On a collar piece it is radial projection in the naturally directed disk coordinate followed by $f$. Thus $d_U$ is directed. Its section $s_U$ is directed because it is $j$ with codomain the induced subspace $U$.

The open sets $V$ and $W$ are disjoint from $A$. The restrictions of $p$ over these sets are therefore bijective quotient maps, hence homeomorphisms. They identify $\abs V$ and $\abs W$ with the products of $\abs T$ with the open disk $\|z\|<b$ and the annulus $a<\|z\|<b$, respectively. In either product $\chi$ sends directed paths to directed paths, and every directed path in the image subdivides into pieces from the cell. These identifications therefore give the product directed structures. In these coordinates $d_W$ and $d_V$ act on the disk coordinate by radial projection and by $\Phi$, and their sections act by multiplication by $c$ and $a$, respectively. All four maps are continuous. They preserve directed paths since they fix the $T$-coordinate and the disk coordinate has its natural directed structure.

\begin{figure}[htbp]
\centering
\begin{tikzpicture}[font=\small,line cap=round,line join=round]
 \begin{scope}
  \fill[orange!16] (0,0) circle (1.22);
  \fill[blue!18] (0,0) circle (.89);
  \fill[blue!5] (0,0) circle (.52);
  \draw (0,0) circle (1.22);
  \draw[dashed] (0,0) circle (.52) (0,0) circle (.89);
  \draw[thin,->] (0,0) -- (1.53,0) node[right] {$\|z\|$};
  \foreach \pos/\val in {.52/a,.89/b,1.22/1}{\draw (\pos,-.05)--(\pos,.05);\node[below,fill=white,inner sep=1pt] at (\pos,-.06) {$\val$};}
  \node at (-.05,.21) {$V$};
  \node at (0,.70) {$W$};
  \node at (0,1.05) {$U$};
 \end{scope}
 \begin{scope}[shift={(3.25,-1.12)}]
  \draw[->] (0,0)--(2.65,0) node[right] {$\|z\|$};
  \draw[->] (0,0)--(0,2.58) node[above] {$\|\Phi(z)\|$};
  \draw[thick,blue!60!black] (0,0)--(1,2.03)--(2.38,2.03);
  \draw[densely dashed] (1,0)--(1,2.03);
  \node[below] at (1,0) {$a$}; \node[below] at (2.38,0) {$1$}; \node[left] at (0,2.03) {$1$};
 \end{scope}
\end{tikzpicture}
\caption{The attaching disk for $m=2$ and the change in radius $\|z\|\mapsto\min(\|z\|/a,1)$. The sets $U$ and $V$ meet the disk in $a<\|z\|\leq1$ and $\|z\|<b$, respectively. The dashed circles bound $W$ and the outer circle attaches to $A\subset U$.}\label{fig:reversible-collar}
\end{figure}
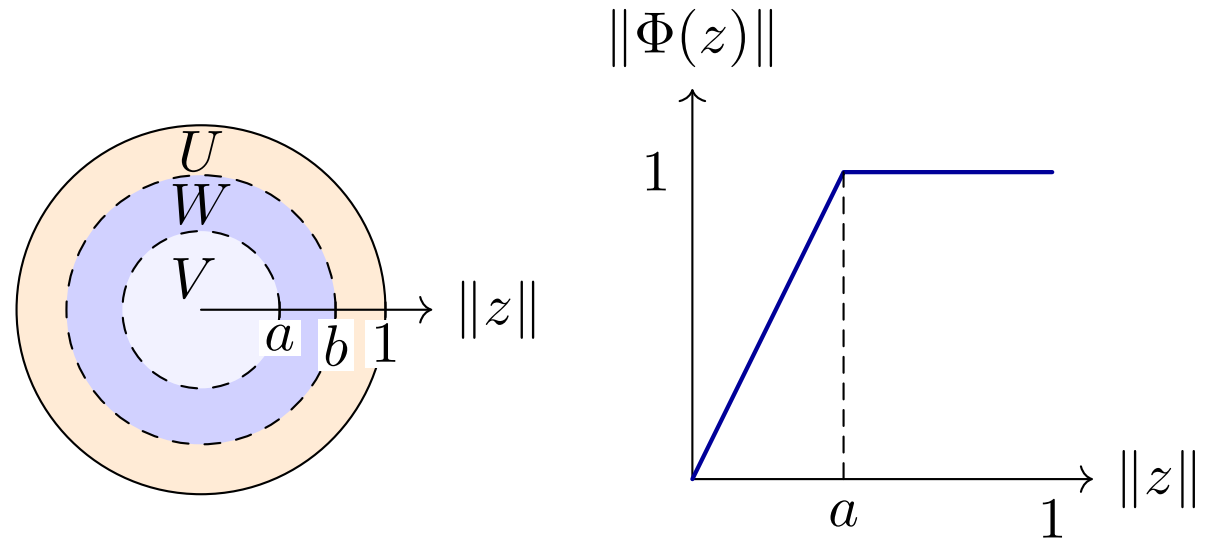

\emph{The homotopies.} On $A\subset U$ put $H_U(x,s)=x$. On the cell pieces of $Z=U,W,V$ put
\[
 H_Z(\chi(t,z),s)=\chi\bigl(t,(1-s)z+s z'\bigr),\qquad
 z'=\begin{cases}
 z/\|z\| & \text{on the collar of }U,\\
 cz/\|z\| & \text{on }W,\\
 a\Phi(z) & \text{on }V .
 \end{cases}
\]
For $s\in[0,1]$ the three segments interpolate radially from $\|z\|$ to $1$, $c$ and $\min(\|z\|,a)$, respectively. They therefore stay in their indicated domains, including at $z=0$ in $V$. In $U$ the endpoint at radius $1$ is identified with $f(t,z/\|z\|)$. Thus for each $Z$ and $x\in Z$,
\[
 H_Z(x,0)=x,\qquad H_Z(x,1)=s_Z(d_Z(x)).
\]

On $W$ and $V$ the interpolations are directed homotopies with parameter $\Iv$ because the disk coordinate carries its natural directed structure. On $U$ the collar interpolation fixes $z$ when $\|z\|=1$, so it agrees with the projection on $A$ at $\chi(t,z)=jf(t,z)$. The two continuous maps thus factor as
\[
 \bigl(\abs A\sqcup\chi^{-1}(\abs U)\bigr)\times I
 \xrightarrow{\ p_U\times\mathrm{id}_I\ }
 \abs U\times I \xrightarrow{\ H_U\ } \abs U .
\]
Since $I$ is locally compact Hausdorff, $p_U\times\mathrm{id}_I$ is a quotient map~\cite[Section~29, Exercise~11(a)]{zbMATH01461253}, so $H_U$ is continuous.\footnote{The products in this factorization are taken in the category of all topological spaces, with the ordinary product topology. The underlying space $\abs{U\times\Iv}$ of the directed product carries the $\Delta$-kelleyfication of $\abs U\times I$. The identity $\abs{U\times\Iv}\to\abs U\times I$ is therefore continuous, and composing it with $H_U$ gives continuity on the directed product.}

For a directed path $(\gamma,\sigma)$ in $U\times\Iv$, subdivide $\gamma$ into directed pieces $j\alpha$ or $\chi(\tau,z)$ as above. On these pieces the image paths are given pointwise by
\[
 H_U(j\alpha,\sigma)=j\alpha,\qquad
 H_U(\chi(\tau,z),\sigma)=\chi\!\left(\tau,(1-\sigma)z+\sigma\frac{z}{\|z\|}\right).
\]
Here $\alpha$ is directed in $A$, $\tau$ is directed in $T$, and the new disk coordinate is continuous, hence directed in $\rev{\mathsf{D}^m}$. Thus each piece is directed in $U$, as is their concatenation.

The homotopy $H_U$ extends to all of $X$ by
\[
 H_X(jx,s)=jx,\qquad H_X(\chi(t,z),s)=\chi\bigl(t,(1-s)z+s\Phi(z)\bigr).
\]
Indeed $\Phi(z)=z/\|z\|$ on the collar and $\Phi(z)=z$ on the attaching sphere. The two continuous formulas therefore descend along the quotient map $p\times\mathrm{id}_I$, as for $H_U$. A directed path in $X\times\Iv$ again subdivides into pieces with $X$-projection in $A$ or from the cell. On each cell piece the $T$-coordinate is preserved and the disk coordinate is continuous. On $A$ the homotopy is the projection. Thus $H_X$ is directed. Set $h(x)=H_X(x,1)$. Then $H_X(x,0)=x$ and
\[
 h|_U=jd_U,\qquad h|_V=\chi d_V .
\]
\end{proof}

\begin{proof}[Proof of Theorem~\ref{thm:disk-excision}]
For $Z=U,W,V$, the section identity and Grandis's invariance theorem applied to $H_Z$ give
\[
 \PiC(d_Z)\PiC(s_Z)=\mathrm{id},\qquad
 \mathrm{id}_{\PiC(Z)}\cong\PiC(s_Z)\PiC(d_Z).
\]
Thus $\PiC(d_Z)$ is an equivalence with quasi-inverse $\PiC(s_Z)$.

Let $P=\PiC(A)\sqcup_{\PiC(B)}\PiC(E)$, and let $q\colon P\to\PiC(X)$ be induced by $\PiC(j)$ and $\PiC(\chi)$. By Theorem~\ref{thm:grandis}, the categorical pushout for the top row of~\eqref{eq:collar-spans} is $\PiC(X)$. In the folk model structure all small categories are cofibrant, and the functors $\PiC(W)\to\PiC(V)$ and $\PiC(B)\to\PiC(E)$ are cofibrations since they are injective on objects. Hovey's cube lemma~\cite[Lemma~5.2.6]{MR99h:55031} therefore implies that the three equivalences $\PiC(d_U)$, $\PiC(d_W)$ and $\PiC(d_V)$ induce an equivalence $d\colon\PiC(X)\to P$ on pushouts.

To identify $q$ itself, recall that $h|_U=jd_U$ and $h|_V=\chi d_V$. The pushout property for the cover therefore gives $qd=\PiC(h)$. Grandis's invariance theorem applied to $H_X$ gives $\PiC(h)\cong\mathrm{id}_{\PiC(X)}$. Since $d$ is an equivalence, so is $q$. The object-set functor preserves colimits, being left adjoint to the functor that gives a set the category with one arrow between each ordered pair of elements. Hence on objects $q$ is the canonical bijection from the pushout of the sets of points to the points of $X$. A fully faithful functor bijective on objects is an isomorphism.
\end{proof}

\begin{proposition}\label{prop:closed-inclusions}
The underlying map of every $\I$-cofibration $i\colon A\to X$ is a closed $T_1$ inclusion, that is, a closed embedding for which every point of $X\setminus i(A)$ is closed in $X$. A path $\gamma$ in $A$ is directed if and only if $i\gamma$ is directed in $X$.
\end{proposition}
\begin{proof}
Each map in~\eqref{eq:generators} is, on underlying spaces, the inclusion of a subcomplex of a finite CW complex. For the strip, $I\times\mathsf{S}^0$ is the union of two opposite edges in the product CW structure on $I\times I$. Thus each underlying map is a $q$-cofibration, since the $q$-cofibrations are generated by the disk boundary inclusions~\cite[Section~2 and Appendix~B, p.~608]{leftproperflow}. The underlying space functor $\D\to\mathbf{Top}$ preserves colimits by Section~\ref{sec:setting} and preserves retracts, so $\abs i$ is a $q$-cofibration. For either choice of $\mathbf{Top}$, disk attachments and their transfinite composites are computed in ordinary topological spaces~\cite[Propositions~B.15 and~B.16]{leftproperflow}. Hence $\abs i$ is a closed $T_1$ inclusion by~\cite[Corollary~2.4.6]{MR99h:55031}.

For the directed structure, first attach a generator $K\hookrightarrow E$ along a directed map $u\colon K\to A$. Write $i\colon A\to X$ and $\chi\colon E\to X$ for the maps to the pushout. In the underlying pushout, $\chi^{-1}(i(A))=K$. Let $\gamma$ be a path in $A$ for which $i\gamma$ is directed. By Section~\ref{sec:setting}, $i\gamma$ subdivides into directed pieces from $A$ or $E$. A piece from $E$ takes values in $K$, which has the directed structure induced by $E$. Applying $u$ gives the corresponding directed piece of $\gamma$ in $A$. Thus $\gamma$ is directed.

For a relative $\I$-cell map, well order the summands at each stage and attach them one at a time. This gives the same pushout by its universal property. In the resulting transfinite sequence $(X_\alpha)$ we prove by induction on $\beta$ that $X_\alpha$ has the directed structure induced by $X_\beta$ for every $\alpha<\beta$. At a successor stage this follows from the single attachment since induced structures restrict successively. At a limit ordinal $\beta$, the underlying space is the union of the earlier stages with its final topology~\cite[Proposition~B.16]{leftproperflow}. The union of their directed path sets is closed under reparametrization and concatenation, since any two stages have a common upper bound below $\beta$. It contains the constant paths, so it is the directed structure of the colimit. A directed path in $X_\beta$ with image in $X_\alpha$ is therefore directed in some $X_\delta$ with $\alpha\leq\delta<\beta$, and hence in $X_\alpha$ by induction.

Finally, the small object argument in the locally presentable category $\D$ exhibits $i\colon A\to X$ as a retract under $A$ of a relative $\I$-cell map $j\colon A\to X'$~\cite[Corollary~2.1.15]{MR99h:55031}. Write $s\colon X\to X'$ and $p\colon X'\to X$ for the maps under $A$, so that $si=j$ and $ps=\mathrm{id}_X$. If $i\gamma$ is directed, then $j\gamma=s\,i\gamma$ is directed in $X'$, so $\gamma$ is directed in $A$. The converse follows because $i$ is a directed map.
\end{proof}

\begin{theorem}\label{thm:excision}
The functor $\PiC$ preserves pushouts along $\I$-cofibrations and colimits of transfinite sequences of $\I$-cofibrations continuous at limit ordinals.
\end{theorem}
\begin{proof}
Let $\lambda$ be a nonzero limit ordinal, let $(X_\alpha)_{\alpha<\lambda}$ be such a sequence, and denote its colimit by $X_\lambda=\operatorname{colim}_{\alpha<\lambda}X_\alpha$. By Proposition~\ref{prop:closed-inclusions}, the transition maps are closed $T_1$ inclusions on underlying spaces, and each stage has the directed structure induced by $X_\lambda$. For either choice of $\mathbf{Top}$ the underlying colimit is the union with its final topology~\cite[Proposition~B.16]{leftproperflow}. Hence every continuous map from the compact interval or square into $\abs{X_\lambda}$ factors through a stage~\cite[Proposition~2.4.2]{MR99h:55031}.

By the induced directed structure, a directed path or elementary $d$-homotopy factors through that stage as a directed map. Thus each arrow of $\PiC(X_\lambda)$ has a representative at one stage. If two representatives give the same arrow, choose a finite zigzag of elementary $d$-homotopies joining them. The maximum of the stage indices for the two paths and the squares of this zigzag is still below $\lambda$. The zigzag in that stage identifies the two classes. On objects the colimit is the union of the point sets. Since filtered colimits of categories are computed on objects and arrows, the canonical functor
\[
 \operatorname{colim}_{\alpha<\lambda}\PiC(X_\alpha)\longrightarrow\PiC(X_\lambda)
\]
is an isomorphism.

By Grandis's invariance theorem, $\PiC$ preserves disjoint unions, including the empty coproduct. Thus $\PiC$ takes the pushout of each generator along any map to a pushout, by the coproduct formula for the generators with empty source and by Theorem~\ref{thm:disk-excision} for the other three.

Let $A\to B$ be a relative $\I$-cell map and let $A\to Y$ be a directed map. Choose an ordinal $\mu$ and a cellular presentation $(B_\alpha)_{\alpha\leq\mu}$ with $B_0=A$ and $B_\mu=B$, attaching one cell at each successor stage. Set $Y_\alpha=Y\sqcup_A B_\alpha$. We prove that the canonical functor $\PiC(Y)\sqcup_{\PiC(A)}\PiC(B_\alpha)\to\PiC(Y_\alpha)$ is an isomorphism. It is the identity on $\PiC(Y)$ at stage zero. If $K\hookrightarrow E$ is attached to $B_\alpha$, its attachment to $Y_\alpha$ along $K\to B_\alpha\to Y_\alpha$ gives $Y_{\alpha+1}$ by the pushout property. Applying the result for a single generator to both attachments gives
\[
\begin{aligned}
 \PiC(Y_{\alpha+1})
 &\cong \PiC(Y_\alpha)\sqcup_{\PiC(K)}\PiC(E)\\
 &\cong \bigl(\PiC(Y)\sqcup_{\PiC(A)}\PiC(B_\alpha)\bigr)\sqcup_{\PiC(K)}\PiC(E)\\
 &\cong \PiC(Y)\sqcup_{\PiC(A)}\PiC(B_{\alpha+1}),
\end{aligned}
\]
where the middle line uses the induction hypothesis and the last uses pasting of pushouts. At a limit ordinal both sequences are colimits of their earlier stages, since pushout along $A\to Y$ preserves colimits under $A$. The successor maps in both sequences are pushouts of generators. Thus $\PiC$ preserves both limit-stage colimits by the first part of the proof. Taking the colimit of the preceding comparison isomorphisms gives the comparison at the limit stage, since pushout along $\PiC(A)\to\PiC(Y)$ also preserves colimits under $\PiC(A)$. For a general $A\to B$ in $\operatorname{cof}(\I)$, take the maps $s\colon B\to B'$ and $p\colon B'\to B$ under $A$ from the proof of Proposition~\ref{prop:closed-inclusions}, where $ps=\mathrm{id}_B$ and $A\to B'$ is a relative $\I$-cell map. By naturality in $B$ these maps exhibit the comparison for $B$ as a retract of the comparison isomorphism for $B'$. Hence it too is an isomorphism.
\end{proof}

\section{The Grandis model structure}\label{sec:model}

For a class $S$ of maps, write $\operatorname{inj}(S)$ for the maps with the right lifting property against $S$, and $\operatorname{cof}(S)$ for its closure under pushouts, transfinite compositions and retracts. The objects of the arrow category are maps, and its morphisms are commutative squares.

\begin{smithRecognition}[{\cite[Definition~1.5 and Theorem~1.7]{MR1780498}}]
Let $\mathcal C$ be a locally presentable category, $I$ a set of maps in $\mathcal C$, and $W$ a class of maps in $\mathcal C$. Assume the following.
\begin{enumerate}
\item[\textup{(c0)}] $W$ satisfies two-out-of-three and is closed under retracts.
\item[\textup{(c1)}] $\operatorname{inj}(I)\subseteq W$.
\item[\textup{(c2)}] $\operatorname{cof}(I)\cap W$ is closed under pushouts and transfinite compositions.
\item[\textup{(c3)}] For each $i\in I$ there is a set $S_i\subseteq W$ such that every morphism $i\to w$ in the arrow category with $w\in W$ factors as $i\to w_0\to w$ for some $w_0\in S_i$.
\end{enumerate}
Then $\mathcal C$ admits a combinatorial model structure with weak equivalences $W$, generating cofibrations $I$, and fibrations $\operatorname{inj}(\operatorname{cof}(I)\cap W)$.
\end{smithRecognition}

Condition \textup{(c3)} is the \emph{solution set condition at $I$}. By~\cite[Definition~1.14 and Proposition~1.15]{MR1780498}, it holds whenever $W$ is an accessible, accessibly embedded full subcategory of the arrow category.

\begin{theorem}[The Grandis model structure]\label{thm:model}
\GrandisModelStatement
\end{theorem}
\begin{proof}
We apply Smith's recognition theorem with $I=\I$ and $W=\WG$. The category $\D$ is locally presentable for either choice of $\mathbf{Top}$ by Section~\ref{sec:setting}. The seven maps of~\eqref{eq:generators} form a set $\I$. Since equivalences of categories are closed under retracts and satisfy two-out-of-three, the same holds for $\WG$. This gives \textup{(c0)}, and Proposition~\ref{prop:lifts} gives \textup{(c1)}.

For $i\in\operatorname{cof}(\I)\cap\WG$, the equivalence $\PiC(i)$ is injective on objects by Proposition~\ref{prop:closed-inclusions}. Thus $\PiC(i)$ is a trivial cofibration in the folk model structure. By Theorem~\ref{thm:excision}, $\PiC$ carries a pushout of $i$ to a pushout of $\PiC(i)$, and a transfinite composite of such maps to the composite of their images. Folk trivial cofibrations are closed under both operations, as is $\operatorname{cof}(\I)$ in $\D$. Hence $\operatorname{cof}(\I)\cap\WG$ is closed under pushouts and transfinite compositions, as required by \textup{(c2)}.

To verify \textup{(c3)} by accessibility, we first show that $\PiC$ is accessible. Since $\D$ is locally presentable, choose a regular cardinal $\kappa$ for which $\cube 0,\cube 1,\cube 2$ are all $\kappa$-presentable. For $i=0,1,2$, let $A_i$ be the functor from $\D$ to sets which sends $X$ to the set $A_i(X)=\D(\cube i,X)$ of directed maps $\cube i\to X$ and acts by postcomposition. Evaluation at the point of $\cube0$ identifies $A_0(X)$ with the set of points of $X$. The map $c_X\colon A_0(X)\to A_1(X)$ sends $x$ to the constant path at $x$. For $j=0,1$ let $r_j\colon A_2(X)\to A_1(X)$ send $F$ to the path $s\mapsto F(s,j)$. We obtain the set $H(X)$ of elementary $d$-homotopies with fixed endpoints as the pullback in sets
\[
\begin{tikzcd}[column sep=large]
 H(X)\ar[r,hook]\ar[d,"{F\mapsto(F(0,0),F(0,1))}"']\ar[dr,phantom,"\lrcorner",very near start]& A_2(X)\ar[d,"{(r_0,r_1)}"]\\
 A_0(X)\times A_0(X)\ar[r,"{c_X\times c_X}"'] & A_1(X)\times A_1(X) .
\end{tikzcd}
\]
For $j=0,1$ let $e_j\colon H(X)\to A_1(X)$ send $F$ to the path $t\mapsto F(j,t)$. For a small category $\mathcal C$, let $\operatorname{Ob}\mathcal C$ and $\operatorname{Mor}\mathcal C$ denote its sets of objects and arrows. With $\operatorname{coeq}(e_0,e_1)$ denoting the coequalizer in sets, we have
\[
 \operatorname{Ob}\PiC(X)=A_0(X),\qquad
 \operatorname{Mor}\PiC(X)=\operatorname{coeq}(e_0,e_1).
\]
Indeed, this coequalizer is the quotient of $A_1(X)$ by the equivalence relation generated by $e_0(F)\sim e_1(F)$ for $F\in H(X)$. Two paths therefore have the same image precisely when a finite zigzag of such homotopies joins them. Each $A_i$ preserves $\kappa$-filtered colimits by $\kappa$-presentability of $\cube i$. Such colimits commute with finite limits in sets and with coequalizers, so the displayed descriptions show that $\PiC$ preserves them on objects and arrows. Since filtered colimits of categories are computed on objects and arrows, $\PiC$ is accessible.

The equivalences of small categories form an accessible, accessibly embedded full subcategory of the arrow category of $\Cat$, by~\cite[Theorem~4.1]{MR2506258} applied to the combinatorial folk model structure. The functor induced by $\PiC$ on arrow categories is accessible because colimits there are computed at source and target. Its inverse image of the equivalences is $\WG$. By~\cite[Definition~1.14 and Proposition~1.18]{MR1780498}, $\WG$ is likewise accessible and accessibly embedded in the arrow category of $\D$.

Thus Beke's Proposition~1.15 gives \textup{(c3)} for $\WG$~\cite{MR1780498}. Smith's theorem now gives the stated model structure.
\end{proof}

\begin{corollary}[Left properness]\label{cor:left-proper}
For either choice of $\mathbf{Top}$, the Grandis model structure $\MG$ is left proper.
\end{corollary}
\begin{proof}
Let $i\colon A\to B$ be a cofibration and $w\colon A\to A'$ a weak equivalence. Form the pushout
\[
\begin{tikzcd}
 A\ar[r,hook,"i"]\ar[d,"w"'] & B\ar[d,"w'"]\\
 A'\ar[r] & B\sqcup_A A'\ar[ul,phantom,"\ulcorner",very near start].
\end{tikzcd}
\]
By Theorem~\ref{thm:excision}, $\PiC(w')$ is the pushout of $\PiC(w)$ along $\PiC(i)$. The functor $\PiC(w)$ is an equivalence of categories, and $\PiC(i)$ is injective on objects by Proposition~\ref{prop:closed-inclusions}, hence a folk cofibration. The folk model structure is left proper since every small category is cofibrant~\cite[Proposition~A.2.4.2]{Lurie2009}. Therefore $\PiC(w')$ is an equivalence, so $w'$ is a weak equivalence in $\MG$.
\end{proof}

\begin{remark}
The domains of the seven generating cofibrations are cofibrant. Besides the initial object they are $\mathsf{S}^0\cong\cube 0\sqcup\cube 0$, $\cube 1\times\mathsf{S}^0\cong\cube 1\sqcup\cube 1$, and $\rev{\mathsf{S}^1}$. The maps from the initial object to $\cube 0$, $\cube 1$ and $\rev{\mathsf{S}^1}$ all occur in~\eqref{eq:generators}. Thus the class of cofibrations is generated by maps between cofibrant objects, as required by Henry. By~\cite[Definition~1.2, Theorem~1.4 and Example~1.7]{henry2020minimal}, there is a minimal model structure with cofibrations $\operatorname{cof}(\I)$. Its weak equivalences belong to $\WG$, so the identity from the minimal to the Grandis model structure is left Quillen. Since the latter has precisely $\WG$ as its weak equivalences, the derived identity localizes the former at the images of all maps in $\WG$. In this sense the Grandis model structure is the left Bousfield localization of the minimal one at $\WG$.

No left Bousfield localization of the discrete model structure on $\D$ has precisely $\WG$ as its weak equivalences. Indeed, every map is a cofibration in such a localization, so its weak equivalences are its trivial cofibrations and are closed under pushout. The collapse $\Iv\to *$ belongs to $\WG$ because $\PiC(\Iv)$ has exactly one arrow between any two points. Let $q\colon\Iv\to\rev{\mathsf{S}^1}$ be given by $q(t)=e^{2\pi it}$. In any cocone on the span below, the map from the circle is constant because its composite with the surjection $q$ is constant. Thus the square is a pushout.
\[
\begin{tikzcd}
 \Iv\ar[r]\ar[d,"q"'] & *\ar[d,"\mathrm{id}"]\\
 \rev{\mathsf{S}^1}\ar[r] & *\ar[ul,phantom,"\ulcorner",very near start].
\end{tikzcd}
\]
The bottom map lies outside $\WG$ because $\PiC(\rev{\mathsf{S}^1})$ is the ordinary fundamental groupoid of the circle~\cite[Section~3.1]{mg} and the collapse sends its vertex group $\mathbb Z$ to the trivial group.
\end{remark}

\begin{remark}\label{rem:interval-boundary}
The boundary inclusion $\mathsf{S}^0\hookrightarrow\cube1$ is not a cofibration in $\MG$. Let $X=[-1,1]$ and put $p(x)=\max(0,x)$. A continuous path $\gamma$ in $X$ is directed when $p\gamma$ is nondecreasing. Thus $p\colon X\to\cube1$ collapses $[-1,0]$ with its natural directed structure. In the following squares $i\in\I$ and $e(0)=-1$, $e(1)=1$.
\[
\begin{tikzcd}
 A\ar[r,"a"]\ar[d,hook,"i"'] & X\ar[d,"p"]\\
 B\ar[r,"b"']\ar[ur,dashed,"\widetilde b"] & \cube1
\end{tikzcd}
\qquad
\begin{tikzcd}
 \mathsf{S}^0\ar[r,"e"]\ar[d,hook] & X\ar[d,"p"]\\
 \cube1\ar[r,"\mathrm{id}"'] & \cube1
\end{tikzcd}
\]
The inclusion $s\colon\cube1\hookrightarrow X$, $s(y)=y$, is a directed section of $p$. When $A=\varnothing$ we can therefore take $\widetilde b=sb$. Otherwise $i$ has the form $T\times\rev{\mathsf{S}^{m-1}}\hookrightarrow T\times\rev{\mathsf{D}^m}$ of Section~\ref{sec:excision}. Any two points of the disk are joined by directed paths in both directions, so the map $z\mapsto b(t,z)$ is constant with value $h(t)=b(t,0)$ for each $t\in T$. The fibre over $0$ is $[-1,0]$, and the fibre over $y\in{]}0,1]$ is $\{y\}$. Both are convex, so for $u\in\mathsf{S}^{m-1}$ we can join $s(h(t))$ to $a(t,u)$ in the fibre over $h(t)$ by
\[
 \widetilde b(t,ru)=(1-r)s(h(t))+r\,a(t,u),\qquad 0\leq r\leq1.
\]
At $r=0$ the value $s(h(t))$ is independent of $u$. The estimate $|\widetilde b(t,ru)-s(h(t))|\leq2r$ gives continuity there, and polar coordinates give continuity for $r>0$. The map $\widetilde b$ extends $a$ and is directed since $p\widetilde b=b$. Thus $p$ is a trivial fibration by Theorem~\ref{thm:model}.

On the right, a lift would send $0$ to $-1$ and every $t>0$ to $t$. The inverse image of the open subset $[-1,-1/2[$ of $X$ would be $\{0\}$, which is not open in $[0,1]$. The directed path $t\mapsto 2t-1$ joins the required endpoints and projects to the identity path with an initial pause, which represents the same arrow in $\PiC(\cube1)$ by~\eqref{eq:thin}.
\end{remark}

\section{A homotopical interpretation of directed Van Kampen}\label{sec:van-kampen}

For a model category $\mathcal M$, write $\mathcal M_\infty$ for its $\infty$-categorical localization at the weak equivalences. By a homotopy colimit in $\mathcal M$ we mean a colimit in $\mathcal M_\infty$. We use this convention for $\MG$ and for $\Cat$ with its folk model structure. For a small category $J$, write $\mathcal M^J$ for the category of functors $J\to\mathcal M$ and natural transformations. For either choice of $\mathcal M$, the localization of $\mathcal M^J$ at the objectwise weak equivalences is equivalent to the functor $\infty$-category $\operatorname{Fun}(NJ,\mathcal M_\infty)$ from the nerve $NJ$ of $J$ to $\mathcal M_\infty$~\cite[Theorem~7.9.8]{Cisinski2019Higher}. Indeed, both model categories admit small colimits, and their cofibrations and trivial cofibrations are closed under coproducts and sequential compositions, as required by the dual of Definition~7.7.2 and by Remark~7.9.7 of~\cite{Cisinski2019Higher}. Under this equivalence the constant diagram functors correspond. Their left adjoints therefore identify the left derived colimit with the colimit in $\mathcal M_\infty$, by the dual of~\cite[Remark~7.9.10]{Cisinski2019Higher}. In particular, every span in $\mathcal M_\infty$ can be represented by an ordinary span in $\mathcal M$.

\begin{theorem}[Directed Van Kampen]\label{thm:positive}
\VanKampenStatement
\end{theorem}
\begin{proof}
By Theorem~\ref{thm:model}, $\PiC$ sends weak equivalences in $\MG$ to equivalences of categories and so induces a functor on the localizations. These localizations admit small colimits~\cite[dual of Proposition~7.7.4]{Cisinski2019Higher}. It suffices to prove that the induced functor preserves pushouts and small coproducts~\cite[Proposition~4.4.2.7]{Lurie2009}.

Consider a span $B\xleftarrow{f}A\xrightarrow{g}C$ in $\MG$. Choose a cofibrant replacement $q_A\colon A'\to A$ by a trivial fibration. Factor $fq_A$ and $gq_A$ as cofibrations $A'\hookrightarrow B'$ and $A'\hookrightarrow C'$ followed by trivial fibrations $q_B\colon B'\to B$ and $q_C\colon C'\to C$. Since $A'$ is cofibrant and both legs are cofibrations, the replacement span is Reedy cofibrant. Hence the pushout $P=B'\sqcup_{A'}C'$ computes the homotopy pushout of the original span~\cite[Theorem~5.1.3 and Corollary~5.1.6]{MR99h:55031}.

By Theorem~\ref{thm:excision}, applying $\PiC$ to the pushout on the left gives the pushout on the right.
\[
\begin{tikzcd}
 A'\ar[r,hook]\ar[d,hook]&B'\ar[d]\\
 C'\ar[r]&P\ar[ul,phantom,"\ulcorner",very near start]
\end{tikzcd}
\qquad
\begin{tikzcd}
 \PiC(A')\ar[r]\ar[d]&\PiC(B')\ar[d]\\
 \PiC(C')\ar[r]&\PiC(P)\ar[ul,phantom,"\ulcorner",very near start].
\end{tikzcd}
\]
The two functors out of $\PiC(A')$ are injective on objects by Proposition~\ref{prop:closed-inclusions}. They are therefore folk cofibrations, and $\PiC(A')$ is cofibrant. The right-hand square computes a folk homotopy pushout. The functors $\PiC(q_A)$, $\PiC(q_B)$ and $\PiC(q_C)$ give an objectwise equivalence from the replacement span to $\PiC(B)\leftarrow\PiC(A)\to\PiC(C)$. Transporting the pushout cocone along this equivalence in the localization exhibits $\PiC(P)$ as the homotopy pushout of the latter span.

For a family $(X_i)_{i\in S}$ indexed by a small set $S$, choose cofibrant replacements $QX_i\to X_i$. Their coproduct computes the homotopy coproduct in $\MG$. By the coproduct formula for $\PiC$ recalled in the proof of Theorem~\ref{thm:excision}, the canonical functor
\[
 \coprod_{i\in S}\PiC(QX_i)\longrightarrow\PiC\!\left(\coprod_{i\in S}QX_i\right)
\]
is an isomorphism, including when $S$ is empty. Each $\PiC(QX_i)$ is folk cofibrant and equivalent to $\PiC(X_i)$, so the left side computes the folk homotopy coproduct.

For an ordinary diagram $X\colon J\to\MG$, the comparison in the statement can now be represented explicitly. Choose projectively cofibrant replacements $QX\to X$ in $\MG^J$ and $v\colon R\to\PiC\circ QX$ in $\Cat^J$. Both projective model structures exist by~\cite[Proposition~A.2.8.2]{Lurie2009}, since the two model categories are combinatorial. The constant diagram functor preserves fibrations and trivial fibrations in these structures, so its left adjoint $\operatorname{colim}_J$ computes the derived colimit on cofibrant diagrams. Applying $\PiC$ to the colimit cocone of $QX$ gives the second functor in
\[
 \operatorname{colim}_J R\xrightarrow{\operatorname{colim}_J v}\operatorname{colim}_J(\PiC\circ QX)
 \longrightarrow\PiC\!\left(\operatorname{colim}_J QX\right).
\]
The composite represents the comparison and is an equivalence of categories by preservation of colimits in the localizations.
\end{proof}

\begin{corollary}[Homotopy pushouts of unions]\label{cor:cover}
Let $U,V$ be induced directed subspaces of $X$ whose sets of points cover $X$. The square of inclusions with vertices $U\cap V,U,V,X$ is a homotopy pushout in $\MG$ if and only if the functor
\[
 \kappa\colon\PiC(U)\sqcup_{\PiC(U\cap V)}\PiC(V)\longrightarrow\PiC(X)
\]
induced by the inclusions into $X$ is an isomorphism. Whenever the interiors of $U$ and $V$ cover $X$, the following square is both a pushout of directed spaces and a homotopy pushout in $\MG$.
\[
\begin{tikzcd}
 U\cap V\ar[r,hook]\ar[d,hook]\ar[dr,phantom,"\lrcorner",very near start]&U\ar[d,hook]\\
 V\ar[r,hook]&X\ar[ul,phantom,"\ulcorner",very near start].
\end{tikzcd}
\]
\end{corollary}
\begin{proof}
Put $W=U\cap V$ with its induced structure, so $W$ is the pullback of $U\hookrightarrow X\hookleftarrow V$. Replace $U\leftarrow W\to V$ by a cofibrant span $U'\leftarrow W'\to V'$ as in the proof of Theorem~\ref{thm:positive}. Its pushout $P$ comes with the comparison $q\colon P\to X$ induced by $U'\to U\hookrightarrow X$ and $V'\to V\hookrightarrow X$. The original square is a homotopy pushout precisely when $q$ is a weak equivalence. By Theorem~\ref{thm:excision}, $\PiC(P)$ is the pushout of the fundamental categories of the replacement span. Both categorical spans are Reedy cofibrant, since both have legs injective on objects. Pushout is left Quillen on spans, so it sends the objectwise equivalence between them to the diagonal equivalence in the commutative triangle
\[
\begin{tikzcd}[column sep=small]
 \PiC(P)\ar[rr,"\PiC(q)"]\ar[dr,"\simeq"']&&\PiC(X)\\
 &\PiC(U)\sqcup_{\PiC(W)}\PiC(V)\ar[ur,"\kappa"']&
\end{tikzcd}.
\]
By Theorem~\ref{thm:model}, $q$ is a weak equivalence precisely when $\PiC(q)$ is an equivalence, hence precisely when $\kappa$ is an equivalence. The object-set functor preserves pushouts, so on objects $\kappa$ is the bijection from the pushout of the sets of points to their union. Thus $\kappa$ is an equivalence precisely when it is an isomorphism.

Suppose now that the interiors cover $X$. Let directed maps $h_U\colon U\to Z$ and $h_V\colon V\to Z$ agree on $W$, and let $h\colon\abs X\to\abs Z$ be the resulting map on points. For every continuous test simplex $\sigma\colon\Delta^n\to\abs X$, a sufficiently fine barycentric subdivision is subordinate to the inverse images of the two interiors, by the Lebesgue number lemma. Each closed simplex of the subdivision is $\Delta$-generated, so the restriction of $\sigma$ factors continuously through the corresponding induced subspace. Pasting on the finite subdivision gives continuity of $h\sigma$, hence of $h$. A directed path in $X$ likewise subdivides into directed paths in $U$ or $V$. Their images under $h_U$ or $h_V$ are directed, so the whole image path is directed. This proves the pushout property in $\D$. Theorem~\ref{thm:grandis} makes $\kappa$ an isomorphism, and the first assertion gives the homotopy pushout.
\end{proof}

\begin{example}[Two closed subintervals]\label{ex:closed-halves}
Fix $0<c<1$ and give $U=[0,c]$ and $V=[c,1]$ their induced directed structures in $\cube1$. Their interiors miss $c$. By~\eqref{eq:thin} and increasing affine rescaling, the fundamental categories of $U$ and $V$ are the corresponding ordered intervals. In their categorical pushout over $\{c\}$, every arrow is a finite composite of arrows from the two intervals. Such a composite is nondecreasing, so if its endpoints lie on one side of $c$, it reduces to the unique arrow in that interval. For $x<c<y$, a composite from $x$ to $y$ must pass through the common object $c$. Composing the factors on either side gives precisely
\[
 x\longrightarrow c\longrightarrow y .
\]
Thus the pushout has exactly one arrow $x\to y$ for $x\leq y$ and none for $x>y$. The comparison with $\PiC(\cube1)$ is an isomorphism, so the union square is a homotopy pushout by Corollary~\ref{cor:cover}.
\end{example}

\begin{example}[Alternating closed subintervals]\label{ex:alternating-cover}
Choose a strictly decreasing sequence $(t_n)_{n\geq0}$ with $t_0=1$ and limit $0$. For $n\geq0$ put $J_n=[t_{n+1},t_n]$. Consider the induced directed subspaces $U,V$ of $\cube1$ given by
\[
 U=\{0\}\cup\bigcup_{k\geq0}J_{2k},\qquad
 V=\{0\}\cup\bigcup_{k\geq0}J_{2k+1}.
\]
Their point sets cover $[0,1]$ and are closed in it, since the intervals are locally finite away from $0$ and both sets contain $0$. Their intersection is $\{0\}\cup\{t_n\mid n\geq1\}$.
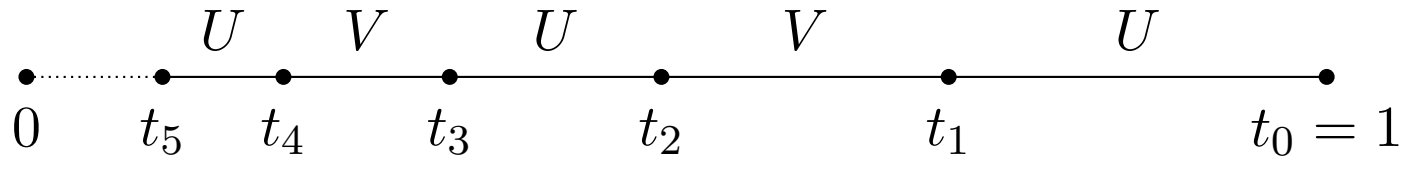
\[
\begin{tikzpicture}[x=1cm,y=1cm,font=\small]
 \draw[densely dotted] (0,0)--(.9,0);
 \draw (.9,0)--(8.6,0);
 \foreach \x/\p in {0/0,.9/{t_5},1.7/{t_4},2.8/{t_3},4.2/{t_2},6.1/{t_1},8.6/{t_0=1}}{
   \fill (\x,0) circle (1.5pt) node[below=2pt] {$\p$};
 }
 \foreach \x/\part in {1.3/U,2.25/V,3.5/U,5.15/V,7.35/U}{
   \node[above=1pt] at (\x,0) {$\part$};
 }
\end{tikzpicture}
\]
For $0<y\leq1$, choose $k$ with $t_{2k}<y$. The midpoints of $J_{2k}$ and $J_{2k+1}$ lie in $]0,y[$ and belong to $U\setminus V$ and $V\setminus U$, respectively. Thus neither piece contains $[0,y]$. A connected subset of either piece containing $0$ is therefore a singleton. In particular every path in either piece starting at $0$ is constant. In fact $0$ is isolated in each $\Delta$-generated subspace topology, since every continuous test simplex whose image contains $0$ has image $\{0\}$.

In $\PiC(U)\sqcup_{\PiC(U\cap V)}\PiC(V)$ an arrow is a finite composite of arrows from the two pieces, so the only arrow with source $0$ is its identity. In $\PiC(\cube1)$ the path $t\mapsto t$ represents an arrow $0\to1$. The comparison functor is therefore not full, and Corollary~\ref{cor:cover} shows that this union square fails to be a homotopy pushout.
\end{example}

\end{document}